\documentclass[10pt]{article}
\usepackage{color,graphicx,wrapfig}

\def\R{\hbox{\bf R}}

\def\N{\hbox{\bf N}}

\def\varphi0{{\cal Q}}

\newcommand{\be}{\begin{equation}}
\newcommand{\ee}{\end{equation}}
\newcommand{\baa}{\begin{array}}
\newcommand{\eaa}{\end{array}}
\newcommand{\ba}{\begin{eqnarray}}
\newcommand{\ea}{\end{eqnarray}}

\newtheorem{theo}{\textsc{Theorem}}[section]

\newtheorem{pro}[theo]{\textsc{Proposition}}
\newtheorem{cor}[theo]{\textsc{Corollary}}
\newtheorem{defi}[theo]{\textsc{Definition}}
\newtheorem{rem}[theo]{\textsc{Remark}}

\begin{document}
\title{{\bf \Large Cross-shaped and Degenerate
Singularities in an\\ Unstable Elliptic Free Boundary Problem}}
\author{
\normalsize\textsc{J. Andersson}\\
{\normalsize\it Max Planck Institute for
Mathematics in the Sciences
  ,}\\
{\normalsize\it  Inselstr. 22, D-04103 Leipzig, Germany}\\
\normalsize\textsc{G.S. Weiss}\\
{\normalsize\it Graduate School of Mathematical Sciences,}\\
{\normalsize\it
University of Tokyo,
3-8-1 Komaba, Meguro, Tokyo,
153-8914 Japan,}\\
{\normalsize\it Guest of the Max Planck Institute
 for Mathematics in the Sciences,}\\
{\normalsize\it Inselstr. 22, D-04103 Leipzig, Germany}
\thanks{G.S. Weiss has been partially supported by the Grant-in-Aid
15740100 of the Japanese Ministry of Education and partially supported
by a fellowship of the Max Planck Society. Both authors thank the Max Planck
Institute for Mathematics in the Sciences for the hospitality
during their stay in Leipzig.}\\
}
%\date{}
\maketitle

%%%%%%%%%%%%%%%%%%%%%%%%%%%%%%%%%%%%%%%%%%%%%%%%%%%%%%%%%%%%%%%%%%%%%%%%%%%%%%%%
%%%%%%%%%%%%%%%%%%%%%%%%%%%%%%%%%%%%%%%%%%%%%%%%%%%%%%%%%%%%%%%%%%%%%%%%%%%%%%%%

\centerline{\small{\bf{Abstract}}}
We investigate singular and degenerate behavior of solutions
of the unstable free boundary problem
$$\Delta u = -\chi_{\{u>0\}}\; .$$
First, we construct a solution that is not of class
$C^{1,1}$ and
whose free boundary consists of four arcs meeting 
in a {\em cross}-shaped singularity.
This solution is completely unstable/repulsive
from above and below which would make
it hard to get by the usual methods, and
even numerics is non-trivial.\\
We also show existence of a degenerate solution. This answers two of the open questions 
in the paper \cite{MW} by R. Monneau-G.S. Weiss.
\noindent{\small{}}\hfill\break

\noindent{\small{\bf{AMS Classification:}}} {\small{35R35, 35J60, 35B65.}}\hfill\break
\noindent{\small{\bf{Keywords:}}} {\small{
free boundary, regularity, cross, asterisk, monotonicity formula, solid combustion, 
singularity, unstable problem}}\hfill\break

%%%%%%%%%%%%%%%%%%%%%%%%%%%%%%%%%%%%%%%%%%%%%%%%%%%%%%%%%%%%%%%%%%%%%%%%%%%%%%%%
%%%%%%%%%%%%%%%%%%%%%%%%%%%%%%%%%%%%%%%%%%%%%%%%%%%%%%%%%%%%%%%%%%%%%%%%%%%%%%%%

\section{Introduction}
We will investigate singular and degenerate
behavior of solutions of the
unstable elliptic
free boundary problem
\begin{equation}\label{eq}
\Delta u = -\chi_{\{ u>0\}}\;\quad  \hbox{ in } \Omega\; .
\end{equation}
The problem (\ref{eq}) is related to traveling wave solutions in solid
combustion with ignition temperature (see the introduction of
\cite{MW} for more details).\\
An equation similar to (\ref{eq}) arises
in the composite membrane problem (see \cite{chanillo1},
\cite{chanillo2}, \cite{blank}).
Another application is the shape of self-gravitating rotating
fluids describing stars (see \cite[equation (1.26)]{stars}). 
\\
This problem has been investigated by R. Monneau-G.S. Weiss in
\cite{MW}. Their main result is that local minimisers of the 
energy
\begin{displaymath}
\int_{\Omega}|\nabla u|^2-2\max(u,0)
\end{displaymath}
are $C^{1,1}$ and that their free boundaries are locally analytic. They also 
establish partial regularity for \textsl{second order non-degenerate}
solutions of (\ref{eq}) (cf. Definition \ref{nondeg}). 
More precisely they show that the
singular set has Hausdorff dimension less than or equal to $n-2$, and that in
two dimensions the free boundary consists close to
singular points of four Lip\-schitz graphs meeting
at right angles. However they left open the question
of the existence of cross-shaped singular points
and of degenerate singularities (cf. \cite[Section 9 and 10]{MW}).

\begin{figure}
\begin{center}
\input{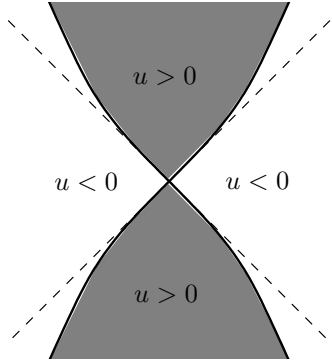}
\end{center}
\caption{A cross-shaped singularity}\label{crossfig}
\end{figure}

\begin{figure}
\begin{center}
\input{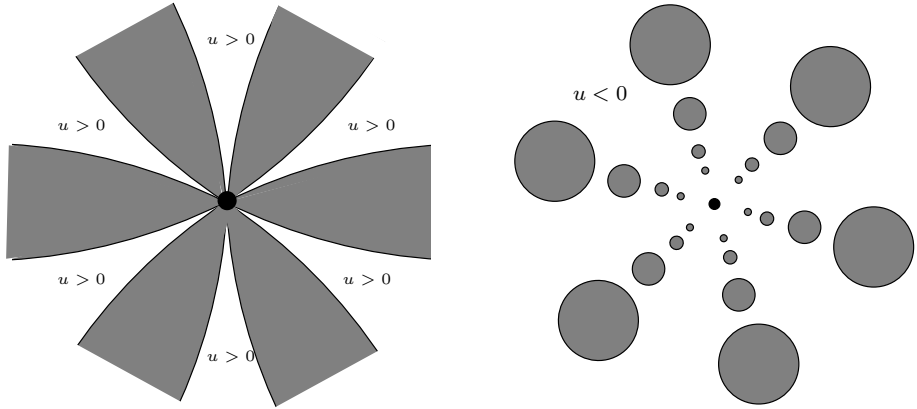}
\end{center}
\caption{Asterisk-shaped singularity or pulse accumulation?}\label{astfig}
\end{figure}

In this paper we will construct both singular points where the free boundary 
consists of four arcs meeting in a cross (see Corollary \ref{cross} and Figure \ref{crossfig})
and solutions that are degenerate of second order at a free boundary point
(see Corollary \ref{ast}).
At this time we do not know whether the shape of the singularity is that
of an asterisk or a product of even higher disconnectivity
(see Figure \ref{astfig}).
\\
In particular, the cross-example
is a counter-example to regularity of the solution
since the solution
is not of class $C^{1,1}$.
\\
In \cite{MW} it has been shown that the second variation
of the energy takes the value $-\infty$ at the function
$x_1^2-x_2^2$. That means that the cross-solution
is completely unstable/repulsive. Moreover it cannot
be approximated from above or below. This makes it
hard to construct it by methods like the implicit
function theorem or comparison methods.\\
Our approach is simple. We construct an operator $T$ such
that each fixed point of $T$,
{\em when adding a certain constant}, satisfies equation (\ref{eq})
{\em and the origin is a point of the $0$-level set!}
By reflection and results from \cite{MW} it is then
possible to show
that origin is non-degenerate of second order
and to obtain the cross.
\\
The construction of degenerate solutions is similar but
simpler.
\\[.4cm]
{\bf Acknowledgement:} We thank Carlos Kenig, Herbert Koch
and R\'egis Monneau for discussions. 
\section{Notation}
Throughout this article $\R^n$ will be equipped with the Euclidean
inner product $x\cdot y$ and the induced norm $\vert x \vert\> .$
We define $e_i$ as the $i$-th unit vector in $\R^n\> ,$ and
$B_r(x^0)$ will denote the open $n$-dimensional ball of center
$x^0\> ,$ radius $r$ and volume $r^n\> \omega_n\> .$
When not specified, $x^0$ is assumed to be $0$.
We shall often use
abbreviations for inverse images like $\{u>0\} := 
\{x\in \Omega\> : \> u(x)>0\}\> , \> \{x_n>0\} := 
\{x \in \R^n \> : \> x_n > 0\}$ etc. 
and occasionally 
we shall employ the decomposition $x=(x_1,\dots,x_n)$ of a vector $x\in \R^n\> .$
When considering a set $A\> ,$ $\chi_A$ shall stand for
the characteristic function of $A\> ,$ 
while
$\nu$ shall typically denote the outward
normal to a given boundary.
\section{Preliminaries}
In this section we state some of the definitions and tools
from \cite{MW}.
\begin{defi}[Non-degeneracy]\label{nondeg}
Let $u$ be a solution of (\ref{eq}) in $\Omega,$
satisfying at $x^0\in \Omega$
\begin{equation}\label{ndeg}
\liminf_{r\to 0} r^{-2}\left(r^{1-n}\int_{\partial B_{r}(x^0)} u^2\> d{\cal H}^{n-1}
\right)^{1\over 2}>0\; .
\end{equation}
Then we call $u$ ``non-degenerate of second order at $x^0$''.
We call $u$ ``non-degenerate of second order'' if it is
non-degenerate of second order at each point in $\Omega$. 
\end{defi}
\begin{rem}
In \cite[Section 3]{MW} it has been shown 
that the maximal solution and each local energy minimiser
are non-degenerate of second order.
\end{rem}
A powerful tool, that we will use in Corollary \ref{cross}, 
is the monotonicity formula introduced 
in \cite{cpde} by one of the authors for a class of semilinear free boundary
problems. For the sake of completeness let us state 
the unstable case here:
\begin{theo}[Monotonicity formula]
\label{mon}
Suppose that $u$ is a solution of (\ref{eq}) in $\Omega$
and that $B_\delta(x^0)\subset \Omega\> .$
Then for all $0<\rho<\sigma<\delta$
the function 
\[ \Phi_{x^0}(r) := r^{-n-2} \int_{B_r(x^0)} \left( 
{\vert \nabla u \vert}^2 \> -\> 2\max(u,0)
\right)\]\[ 
- \; 2 \> r^{-n-3}\>  \int_{\partial B_r(x^0)}
u^2 \> d{\cal H}^{n-1}\; ,\]
defined in $(0,\delta)\> ,$ satisfies the monotonicity formula
\[ \Phi_{x^0}(\sigma)\> -\> \Phi_{x^0}(\rho) \; = \;
\int_\rho^\sigma r^{-n-2}\;
\int_{\partial B_r(x^0)} 2 \left(\nabla u \cdot \nu - 2 \>
{u \over r}\right)^2 \; d{\cal H}^{n-1} \> dr \; \ge 0 \; \; .\] 
\end{theo}
The following proposition has been proven in \cite[Section 5]{MW}.
\begin{pro}[Classification of blow-up limits with fixed center]\label{fixedcenter}
Let $u$ be a solution of (\ref{eq}) in $\Omega$
and let us consider a point
$x^0\in \Omega\cap \{ u=0\}\cap\{ \nabla u =0\}.$\\
1) In the case $\Phi_{x^0}(0+)=-\infty$,
$\lim_{r\to 0} r^{-3-n}\int_{\partial B_r(x^0)} u^2 \> d{\cal H}^{n-1}
= +\infty$, and for
$S(x^0,r) := \left(r^{1-n}\int_{\partial B_{r}(x^0)} u^2\> d{\cal H}^{n-1}
\right)^{1\over 2}$
each limit of
\[ \frac{u(x^0+r x)}{S(x^0,r)}\]
as $r\to 0$ is a homogeneous harmonic polynomial of degree $2$.
\\
2) In the case $\Phi_{x^0}(0+)\in (-\infty,0)$,
\[ u_r(x) := \frac{u(x^0+r x)}{r^2}\]
is bounded in $W^{1,2}(B_1(0))$,
and each limit as $r\to 0$ is a homogeneous solution of degree $2$.\\
3) Else $\Phi_{x^0}(0+)=0$, and
\[ \frac{u(x^0+r x)}{r^2}\to 0\hbox{ in } W^{1,2}(B_1(0)) \hbox{ as } r\to 0\; .\]
\end{pro}
\begin{rem}\label{int}
1) As shown in \cite[Lemma 5.2]{MW}, 
the case 2) is not possible in two dimensions.\\
2) Case 3) is equivalent to $u$ being degenerate of second order
at $x^0$.
\end{rem}
\section{Main Results}
Let $\pi/{\phi_0}\in \N$ and let us define the disk sector $K=K_{\phi_0}=\{r(\cos \phi,\sin\phi):\, 0<r<1, 0<\phi<
\phi_0\}$.
For $g\in C^{\alpha}(\partial B_1\cap \partial K)$, $C^{\alpha}_g(\bar K)$ will denote the 
subspace of $C^{\alpha}(\bar K)$ consisting of 
all the functions with boundary values $g$ on $\partial B_1\cap \partial K$.
\\
Consider now the operator 
$T=T_{\epsilon,g}:C^{\alpha}_g(\bar K)\rightarrow C^{\alpha}_g(\bar K)$ defined by
\begin{displaymath}
\begin{array}{ll}
\Delta T(u)=-f_{\epsilon}(u-u(0)) & \textrm{in } K \> ,\\
T(u)=g & \textrm{on }  \partial B_1\cap \partial K\> ,\\
\frac{\partial (T(u))}{\partial \nu}=0 & \textrm{on } \partial K-\partial B_1\> ;
\end{array}
\end{displaymath}
here $f_\epsilon\in C^\infty(\R),f_\epsilon(z)\ge 
\chi_{\{ z>0\}}$ in $\R$ and $f_\epsilon
\downarrow \chi_{\{ z>0\}}$ as $\epsilon\downarrow 0.$

Since there exists for $F\in L^\infty(K)$
a $W^{1,2}(K)$-solution $v$ of
\begin{displaymath}
\begin{array}{ll}
\Delta v = F & \hbox{in } K\; ,\\
v=g & \textrm{on }  \partial B_1\cap \partial K\> ,\\
\frac{\partial v}{\partial \nu}=0 & \textrm{on } \partial K-\partial B_1\> ,
\end{array}
\end{displaymath}
we obtain after reflection a $W^{1,2}(B_1)$-function
that solves $\Delta v = F$ in $B_1-\{ 0\}$,
where $F$ means the reflected function defined on $B_1$.
As the origin is a set of vanishing capacity, 
$v$ is a weak solution of $\Delta v = F$ in $B_1$.
Applying
the regularity theory for elliptic equations (see 
for example \cite[Lemma 9.29]{GT}), we see that
$T$ is for small $\alpha$ a continuous compact operator from 
$C^{\alpha}_g(\bar K)$ into itself, and that
\begin{displaymath}
\|T_{\epsilon,g}(u)\|_{C^{\alpha}(\bar K)}\le C\; ,
\end{displaymath}
where $C$ is a constant depending only on $g$.

From Schauder's fixed point theorem (see for example \cite[Chapter 11]{GT}) 
we infer that $T_{\epsilon,g}$ has a 
fixed point $u_{\epsilon}\in 
C^{\alpha}_g(\bar K)\cap \{\|u\|_{C^{\alpha}(\bar K)}\le C\}$.
Alternatively, we could also show existence of a fixed
point in a class of symmetric functions.\\
Reflecting and applying
$L^p$-estimates we obtain a sequence $\epsilon_m\to 0$
such that the reflected $u_{\epsilon_m}-u_{\epsilon_m}(0)\to u$ strongly 
in $C^{1,\beta}(\overline{B_{1-\delta}})$ 
and weakly in $W^{2,p}(B_{1-\delta})$
for each $\delta \in (0,1)$
as $m\to \infty$. At a.e. point of $\{ u>0\} \cup \{ u<0\}$,
$u$ satisfies the equation $\Delta u= -\chi_{\{ u>0\}}$.
At a.e. point of $\{ u=0\}$, the weak second derivatives
of the $W^{2,2}$-function $u$ are $0$, so that we obtain:
\begin{pro}[Existence of a Fixed Point]\label{pro:pro}
For each $g\in C^{\alpha}(\partial B_1\cap\partial K)$ there exists a constant
$\kappa$ such
that the boundary value problem
\begin{displaymath}
\begin{array}{ll}
\Delta u=-\chi_{\{u>0\}}& \textrm{in } K \\
u=g-\kappa & \textrm{on } \partial B_1\cap \partial K\> ,\\ 
\frac{\partial u}{\partial \nu}=0 & 
\textrm{on } \partial K-\partial B_1
\end{array}
\end{displaymath}
has a solution $u\in \bigcap_{\delta \in (0,1)}C^{1,\beta}(\bar K\cap \overline{B_{1-\delta}})$ such that $u(0)=0.$
\end{pro}
We will use Proposition \ref{pro:pro}
to prove the existence of singular and degenerate solutions:
\begin{cor}[Construction of a Cross-shaped Singularity]\label{cross}
There exists a solution $u$ of 
\begin{displaymath}
\Delta u= -\chi_{\{u>0\}} \quad \textrm{in } B_1
\end{displaymath}
that is not of class $C^{1,1}$,
such that each limit of
\[ \frac{u(r x)}{S(0,r)}\]
as $r\to 0$ is after rotation the function
$(x_1^2-x_2^2)/\Vert x_1^2-x_2^2\Vert_{L^2(\partial B_1(0))}$.
\end{cor}
\textsl{Proof:} By Proposition \ref{pro:pro} there exists
for each $M\in \R-\{ 0\}$ a constant $\kappa \in \R$ and
a solution in $K_{\pi/2}$ with boundary values $g=M(x_1^2-x_2^2)-\kappa$ 
on $\partial B_1\cap \partial K_{\pi/2}$ 
satisfying the homogeneous Neumann boundary condition on 
$\partial K_{\pi/2}-\partial B_1$.
Using
the homogeneous Neumann boundary condition 
and the fact that $u\in C^{1,\beta}(\overline{K_{\pi/2}\cap B_{1-\delta}})$
we can reflect 
this solution twice at the coordinate axes to obtain 
a solution in the unit ball $B_1$, called again $u$.

Also by Proposition \ref{pro:pro}, we know that 
$u(0)=0$. Thus $u(0)=0$ and $\nabla u(0)=0$ so that
Proposition \ref{fixedcenter} applies.
What remains to be done is to exclude case 3) of
Proposition \ref{fixedcenter} (see Remark \ref{int} 1)). That done, it follows
from the statement in case 1) that $u$ is not
of class $C^{1,1}$.\\
To this end we use the monotonicity formula Theorem \ref{mon}. 
If
$\lim_{r\rightarrow 0}\Phi_0(r)=0$, then
$\Phi_0(r)\ge 0$ for all $r>0$.
Therefore we only need to show that $\Phi_0(1)< 0$:

For $h=M(x_1^2-x_2^2)$ and $g=h$ let us write $u=v+h-\kappa$:
The function $v$ satisfies
\begin{displaymath}
\begin{array}{ll}
\Delta v = \Delta u & \textrm{in } B_1 \hbox{ and}\\
v=0 & \textrm{on } \partial B_1.
\end{array}
\end{displaymath}
Notice that $-1\le \Delta v \le 0$ implies that $0<v<C_1$ and $|\nabla v|<C_1$ 
where $C_1$ is a universal
constant. In particular $C_1$ is independent of $M$. 
We also know that
$\kappa=v(0)\in (0,C_1)$ since $u(0)=0$. 
Now we calculate the energy $\Phi_0(1)$
of $u$.
\begin{displaymath}
\begin{array}{l}
\Phi_0(1)=\int_{B_1}|\nabla u|^2 - 2u^+-2\int_{\partial B_1}u^2\> d{\cal H}^{n-1} \\
= \int_{B_1}|\nabla (v+h)|^2 - 2(v+h-\kappa)^+-2\int_{\partial B_1}
(v+h-\kappa)^2\> d{\cal H}^{n-1} \\
= \int_{B_1}|\nabla v|^2+2\nabla v \cdot \nabla h +|\nabla h|^2 - 2(v+h-\kappa)^+-2\int_{\partial B_1}(h-\kappa)^2\> d{\cal H}^{n-1}\; ,
\end{array}
\end{displaymath}
where we have used that $\kappa$ is a constant and that $v=0$ on $\partial B_1$.
Integrating by parts and using the specific form of $h$ shows that
\begin{displaymath}
\begin{array}{l}
\Phi_0(1)=\int_{B_1}|\nabla v|^2-2(v+h-\kappa)^+-2\int_{\partial B_1}\kappa^2\> d{\cal H}^{n-1} \\
< \int_{B_1}|\nabla v|^2- 2(v+h-\kappa)^+< \int_{B_1}C_1^2-2(h-C_1)^+  \\
=\int_{B_1}C_1^2-2(M(x_1^2-x_2^2)-C_1)^+.
\end{array}
\end{displaymath}
The last integral is negative if $M$ is large. We have thus shown that
$\Phi_0(1)<0$ for sufficiently large $M$. 
\begin{rem}
To calculate the just obtained solution numerically would
-- because of the severe instability -- not be easy.
\end{rem}
The next corollary establishes the existence of degenerate solutions
of second order:
\begin{cor}[Construction of a Degenerate Point]\label{ast}
There exists
a non-trivial solution $u$ of
\begin{displaymath}
\Delta u= -\chi_{\{u>0\}} \quad \textrm{in } B_1
\end{displaymath}
that is degenerate
of second order at the origin.
\end{cor}
\textsl{Proof:} This is also a direct consequence of Proposition \ref{pro:pro}.
The proposition yields a solution in $K_{\pi/4}$ with 
boundary data $\cos(4\phi)-\kappa$ on $\partial K_{\pi/4}\cap \partial B_1$.
Let us reflect this solution three times to get a solution $u$ 
in the unit ball $B_1$.
As in the previous corollary $0=u(0)=|\nabla u(0)|$.
We only have 
to show that $u$ is degenerate of second order. 
Suppose towards a contradiction that this is not true:
then by Remark \ref{int} 1), case 1) of Proposition \ref{fixedcenter}
has to apply. We obtain after a rotation
a blow-up limit of the form
$(x_1^2-x_2^2)/\Vert x_1^2-x_2^2\Vert_{L^2(\partial B_1(0))}$.
But there is no rotation for which
that blow-up limit could be symmetric with respect
to the two axes $x_1=0$ and $x_1=x_2$, yielding a contradiction.
%%%%%%%%%%%%%%%%%%%%%%%%%%%%%%%%%%%%%%%%%%%%%%%%%%%%%%%%%
\section{Open Questions}
Concerning the set of degenerate singular points there
remains the question whether {\em large} degenerate singular
sets are possible. Also it would be nice to know
the precise shape of isolated degenerate singularities,
and whether infinite order vanishing is possible or not.
%%%%%%%%%%%%%%%%%%%%%%%%%%%%%%%%%%%%%%%%%%%%%%%%%%%%%%%%%%%%%%%%%%%%%%%%%%%%%%%%

\bibliographystyle{plain}
\bibliography{anderssonweiss050810.bib}

\def\cprime{$'$} \def\cprime{$'$}
\begin{thebibliography}{1}

\bibitem{blank}
Ivan Blank.
\newblock Eliminating mixed asymptotics in obstacle type free boundary
  problems.
\newblock {\em Comm. Partial Differential Equations}, 29(7-8):1167--1186, 2004.

\bibitem{stars}
Luis~A. Caffarelli and Avner Friedman.
\newblock The shape of axisymmetric rotating fluid.
\newblock {\em J. Funct. Anal.}, 35(1):109--142, 1980.

\bibitem{chanillo2}
S.~Chanillo, D.~Grieser, M.~Imai, K.~Kurata, and I.~Ohnishi.
\newblock Symmetry breaking and other phenomena in the optimization of
  eigenvalues for composite membranes.
\newblock {\em Comm. Math. Phys.}, 214(2):315--337, 2000.

\bibitem{chanillo1}
S.~Chanillo, D.~Grieser, and K.~Kurata.
\newblock The free boundary problem in the optimization of composite membranes.
\newblock In {\em Differential geometric methods in the control of partial
  differential equations (Boulder, CO, 1999)}, volume 268 of {\em Contemp.
  Math.}, pages 61--81. Amer. Math. Soc., Providence, RI, 2000.

\bibitem{GT}
David Gilbarg and Neil~S. Trudinger.
\newblock {\em Elliptic partial differential equations of second order}, volume
  224 of {\em Grundlehren der Mathematischen Wissenschaften [Fundamental
  Principles of Mathematical Sciences]}.
\newblock Springer-Verlag, Berlin, 1983.

\bibitem{MW}
R\'egis Monneau and G.S. Weiss.
\newblock An unstable elliptic free boundary problem arising in solid
  combustion.
\newblock {\em http://arxiv.org/abs/math.AP/0507315, submitted}.

\bibitem{cpde}
Georg~S. Weiss.
\newblock Partial regularity for weak solutions of an elliptic free boundary
  problem.
\newblock {\em Comm. Partial Differential Equations}, 23(3-4):439--455, 1998.

\end{thebibliography}

\end{document}